\documentclass[10pt]{amsart}

\usepackage{amsmath,amsthm,amsfonts,amssymb}
\usepackage{graphicx}
\usepackage[usenames]{color}
\usepackage{array}
\usepackage{colortbl}

\newtheorem{theorem}{Theorem}[section]
\newtheorem{lemma}[theorem]{Lemma}
\newtheorem{proposition}[theorem]{Proposition}
\newtheorem{corollary}[theorem]{Corollary}

\theoremstyle{definition}
\newtheorem{algorithm}[theorem]{Algorithm}

\theoremstyle{remark}
\newtheorem{remark}[theorem]{Remark}

\newcommand{\gp}[1]{{\left\langle #1 \right\rangle}}
\newcommand{\rb}[1]{{\left( #1 \right)}}

\newcommand{\url}[1]{#1}

\def\CT{{\mathcal T}}
\def\CR{{\mathcal R}}

\def\MZ{{\mathbb{Z}}}

\let\al\alpha

\mathcode`\*="0203

\title[The Conjugacy Problem in the Grigorchuk Group]{The Conjugacy Problem in the Grigorchuk Group is polynomial time decidable}

\date{July 20, 2008}

\author[] {Igor Lysenok}
\address{Steklov Mathematical Institute, Moscow, Russia}
\email{igor.lysenok@gmail.com}

\author[]{Alexei Myasnikov}
\address{Department of Mathematics, McGill University,
Montreal, Canada} \email{amiasnikov@gmail.com}

\author[]{Alexander Ushakov}
\address{Department of Mathematics, Stevens Institute of Technology,
Hoboken, USA} \email{sasha.ushakov@gmail.com}

\begin{document}

\maketitle

\pagestyle{myheadings} \markright{{\sl LMU $\bullet$ SCP in Grigorchuk Group,
 $\bullet$ 07.20.08}}

\begin{abstract}
In this paper we prove that the
Conjugacy Problem in the Grigorchuk group $\Gamma$ has polynomial time complexity. This solves  the problem posed in \cite{Grig2006} rather unexpectedly.
\end{abstract}

\tableofcontents

\section{Introduction}
\label{se:Introduction}

In this paper we discuss algorithmic complexity of the conjugacy problem in the original
Grigorchuk group  $\Gamma$. The  group $\Gamma$ first appeared in \cite{Grig80} almost 30 years ago, now it  plays an important part in several areas of modern group theory:
 growth in groups \cite{Grig91}, Burnside's problems \cite{Grig80},   amenability \cite{Grig96}, just infinite
groups \cite{Grig99}.
Recently the group $\Gamma$ was proposed as a possible
platform for cryptographic schemes (see \cite{GZ,P,MU5}), where the algorithmic security
of the schemes is based on  the computational
hardness of certain variations of the word and conjugacy problems in $\Gamma$.
Bibliography on $\Gamma$ is quite extensive, here we refer to
publications \cite{Harpe,Grig2006} that give a  comprehensive and
accessible survey on $\Gamma$.

Our interest in $\Gamma$ comes from rather different direction, it concerns with foundations of algorithmic group theory.   Recall that  the classical approach to algorithmic problems in groups deals mostly with finitely
  presented groups --  an old tradition, coming from topology.   Another way to study  algorithmic problems in groups stems from constructive mathematics, where elements of a group have  to be given as finitary  objects (matrices over number fields, automorphisms of graphs, complexes, or other constructible objects) and the group multiplication has to be effectively described or computable -- Rabin's recursive groups \cite{Rabin}  or Malcev's constructible groups \cite{Malcev} provide typical examples of this type. A more general approach to algorithmic group theory concerns with groups given by arbitrary  recursive presentations. There are some known general results in this direction, including the spectacular Higman's embedding theorem \cite{Higman}, but a cohesive  theory is lacking (perhaps, due to the huge variety of groups in this class).  The Grigorchuk group $\Gamma$ may serve as a model case of study. Indeed,  $\Gamma$   can be easily  described as generated by  four particular automorphisms of the infinite rooted binary tree,  but  it is not finitely presented, though it has a nice infinite recursive ``self-similar''  presentation.
  Studying algorithmic problems in $\Gamma$  may  provide some interesting insights on how to deal with
  recursively presented groups whose  presentations  are infinite
  but can be described by repeating some typical finite pattern or
  obvious self-similarity.

The Word, Conjugacy, and Isomorphism are the  three famous Dehn's  algorithmic problems in group theory.   The Word Problem in $\Gamma$ is decidable and its time complexity is  $O(n \log n)$ (see, for example, \cite{Harpe,Grig2006}).  It has been  shown in \cite{Leonov,Rozjkov} that the Conjugacy Problem (CP) is decidable in $\Gamma$.
In fact, $\Gamma$ is conjugacy separable \cite{WZ}. Moreover, \cite{Leonov} gives  a complete characterization of Grigorchuk groups $G_\omega$ with decidable CP -- precisely those ones where the  sequence $\omega$ is recursive.
Another decision algorithm for CP in $\Gamma$ is described in \cite{BGS03} and
\cite{Grig2006}.  This is a {\em branching} algorithm, it  is based on a
branching rewriting process, similar to the original decision algorithm
for the Word Problem in $\Gamma$  \cite{Grig80}. The time upper bound for
this algorithm given in \cite{Grig2006} is double exponential.
This raises a natural question (see Problem 5.1 in  \cite{Grig2006}): what is the time complexity
of CP in $\Gamma$?

We show below
that CP in $\Gamma$ can be solved in polynomial time. To prove this we modify
the decision algorithm from  \cite{Grig2006}:  given two elements
$u, v \in \Gamma$ we construct, first,  a unique conjugacy tree $T_{u,v}$ (there were
exponentially many trees in \cite{Grig2006}), then we provide a
routine, similar to the one in \cite{Grig2006}, which given a
conjugacy tree $T_{u,v}$ checks whether $u$ and $v$ are
conjugate in $G$ or not. This routine requires polynomial
time in the size of $T_{u,v}$. Finally, we show, and this the main technical result of the paper,
 that the size
of $T_{u,v}$ is polynomial in the total length $|u| + |v|$, so
the decision  algorithm is polynomial in time.
This part is tricky, to prove it  we replace, following \cite{Bartholdi:1998},  the standard length on $\Gamma$ by a new,  "weighted" length, called the {\em norm}, and show that the standard splitting $w \rightarrow (w_0,w_1)$ of elements from $St_\Gamma(1)$  has very nice metric properties relative to the norm.
These  metric properties allow one to prove that the  length of the elements that appear in the construction of $T_{u,v}$ drops exponentially, so the height of the tree $T_{u,v}$ is about $log(|u| +|v|)$, hence the size of $T_{u,v}$ is polynomial in $|u| +|v|$. The degree $d$ of the polynomial depends on the metric properties of the splitting,  currently $d = 7$.
The resulting decision algorithm for CP in $\Gamma$ has the upper time bound $O(n^{8})$. We would like to point out that it is not clear whether this  upper bound is tight or not. In fact, all our computer experiments indicate that the algorithm is quite practical, it behaves like an algorithm with a quadratic time upper bound. The algorithm itself is available on line \cite{CRAG}.
 Finally, we want to mention that it seems plausible  that a similar method could give a polynomial time decision algorithm for CP in some other self-similar contracting groups.

\section{Preliminaries on the Grigorchuk group}
\label{se:preliminaries}

In this section, following \cite{Harpe} and \cite{Grig2006}, we define the
Grigorchuk group $\Gamma$ and recall some of its properties. Notation and the techniques introduced here are heavily used throughout the paper.

\subsection{The Grigorchuk group $\Gamma$}

For a set $X$ by $X^\ast$ we denote the set of all finite words (sequences)  in $X$. If $u \in X^\ast$ and $x \in X$ then $|u|_x$ is the {\em number of occurrences} of $x$ in $u$ and $|u|$ is the {\em length} of $u$.

\begin{figure}[h]
\centerline{\includegraphics[scale=0.6]{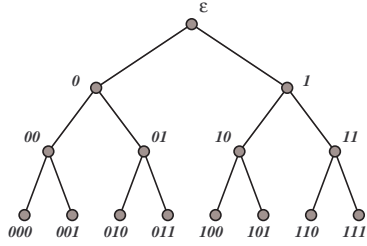} }
\caption{Binary tree with labeled vertices.}
\label{fi:bin_tree}
\end{figure}

Let $\CT$ be an infinite  rooted regular binary tree. Recall that the vertex set of $\CT$ is precisely the set $\{0,1\}^\ast$ of all finite binary words (the empty word  $\varepsilon$ at the root) and two words $u$ and $v$ are connected by an edge in $\CT$ if and only if one of them, say $v$, is obtained from the other by adding one bit $b \in \{0,1\}$ at the end, so $v=ub$. The tree $\CT$ is shown in Figure \ref{fi:bin_tree}.

Let $Aut(\CT)$ be the group of automorphisms of $\CT$ as a rooted tree. Note that any automorphism
of $\CT$ fixes the root $\varepsilon$. Clearly, every $\phi \in Aut(\CT)$ either fixes vertices $0, 1$ or permutes them. The ones that fix $0$ and $1$ form a normal subgroup $St(1)$ of $Aut(\CT)$ of index $2$. Let $T_0$ be the ``left'' subtree of $\CT$, i.e., the subgraph induced by all vertices that start with $0$, and $T_1$ the ``right'' subtree of $\CT$ induced by all vertices starting with $1$. The  automorphism $a \in Aut(\CT)$, defined on vertices of $\CT$ by
    $$a(b_1, b_2, \ldots, b_n) = (1-b_1, b_2, \ldots, b_n),$$
 swaps  the subtrees $T_0$ and $T_1$, hence $a \not \in St(1)$ and $Aut(\CT) = St(1) \sqcup St(1)a$.

The Grigorchuk group $\Gamma$ is the  subgroup of $Aut(\CT)$ generated by four automorphisms
$a,b,c,d$, where $b,c, d \in St(1)$ are defined recursively as follows:
$$
{\bf b}(b_1, b_2, \ldots, b_n) =
\left\{
\begin{array}{ll}
(b_1, 1-b_2, b_3,\ldots, b_n),  & \mbox{if } b_1=0; \\
(b_1, {\bf c}(b_2, \ldots, b_n)), & \mbox{if } b_1=1; \\
\end{array}
\right.
$$
$$
{\bf c}(b_1, b_2, \ldots, b_n) =
\left\{
\begin{array}{ll}
(b_1, 1-b_2, b_3,\ldots, b_n),  & \mbox{if } b_1=0; \\
(b_1, {\bf d}(b_2, \ldots, b_n)), & \mbox{if } b_1=1; \\
\end{array}
\right.
$$
$$
{\bf d}(b_1, b_2, \ldots, b_n) =
\left\{
\begin{array}{ll}
(b_1, b_2, \ldots, b_n),    & \mbox{if } b_1=0; \\
(b_1, {\bf b}(b_2, \ldots, b_n)), & \mbox{if } b_1=1; \\
\end{array}
\right.
$$

It is easy to see that the automorphisms $a,b,c,d$ satisfy the relations
\begin{equation}\label{le:grig_relations}
    a^2=b^2=c^2=d^2=1, ~ bc=cb=d.
\end{equation}
In particular,
    $$\gp{a} \simeq \MZ_2 ~~\mbox{ and }~~ \langle b,c,d \rangle \simeq \MZ_2 \times \MZ_2.$$
Consider a group
    $$\Gamma_0 = \gp{a,b,c,d \mid a^2=b^2=c^2=d^2=1, ~ bc=cb=d}.$$
The group $\Gamma_0$ is the  free product of the subgroups $\langle a\rangle$ and $\langle b,c,d \rangle$.
It follows  that any word
$w \in (X \cup X)^\ast$ is equal in $\Gamma_0$ to a unique {\em reduced} word
\begin{equation}\label{eq:grig_reducedword}
  red(w) =    u_0 a u_1 \ldots u_n a u_{n+1}
\end{equation}
where $u_1, \ldots, u_n \in \{b,c,d\}$, $u_0, u_{n+1} \in \{\varepsilon, b, c, d\}$, in particular, $u_0, u_{n+1}$ could be empty. The following rewriting system $\mathcal{W}$:
    $$x^2 \rightarrow \varepsilon, \ \  x^{-1} \rightarrow x,   \ \ \ (x \in X)$$
    $$rs \rightarrow t, \ \ \ (r,s,t  \in \{b,c,d\}, r \neq s \neq t)$$
is terminating and confluent, and $red(w)$ is precisely the reduced form of $w$ relative to  $\mathcal{W}$. Clearly, $|red(w)| \leq |w|$. Furthermore, given a word $w \in (X \cup X)^\ast$ one can compute  $red(w)$ in time $O(|w|)$.

Denote by $\CR$ the set of all reduced words in $X^\ast$ and by
$\CR^e$ the set of all reduced words $w$ in $X^\ast$ such that $|w|_a$ is even.

Let $St_\Gamma(1) = St(1) \cap \Gamma$  be the set of automorphisms in $\Gamma$ stabilizing the
first level of $\CT$, i.e., stabilizing the vertices $\{0,1\}$.

\begin{lemma} \label{le:st_generators}
The following holds:
 \begin{enumerate}
 \item [1)] A word $w \in X^\ast$ represents an element of $St_\Gamma(1)$
if and only if $|w|_a$  is even.
 \item [2)] $\Gamma = St_\Gamma(1) \sqcup a St_\Gamma(1)$.
 \item [3)] $St_\Gamma(1) = \langle b,c,d,aba,aca,ada\rangle$.
\end{enumerate}
\end{lemma}

\begin{proof}
Follows from the definition of the group $\Gamma$ and the definition of the
elements $a$, $b$, $c$, and $d$.

\end{proof}

Every automorphism $g \in St(1)$ fixes the first level of $\CT$ and hence induces automorphisms $g_0$ and $g_1$ on the subtrees $T_0$ and $T_1$ of $\CT$. Since the subtrees $T_0$ and $T_1$ are naturally isomorphic to $\CT$  the mapping $g \mapsto (g_0,g_1)$ gives a
group isomorphism
    $$\psi: St(1) \rightarrow Aut(\CT) \times Aut(\CT)$$
and, hence, in the event $g \stackrel{\psi}{\mapsto} (g_0,g_1)$ we can write $g = (g_0,g_1)$.

If $g,h \in St(1)$ and $g = (g_0,g_1)$ and $h = (h_0,h_1)$
then (since $\psi$ is an isomorphism)
\begin{equation}\label{eq:st_split_prod}
    gh = (g_0h_0,g_1h_1)
\end{equation}
and (an easy  computation)
\begin{equation}\label{eq:st_split_conj}
    a^{-1}ga = (g_1,g_0).
\end{equation}

We use these formulas frequently and without references.
Also, it is easy to check that for generators $b,c,d$ the following equalities hold
\begin{equation}\label{eq:st_split}
\left\{
\begin{array}{llllll}
b = (a,c);&&&&&\\
c = (a,d);\\
d = (1,b);\\
\end{array}
\right.
\left\{
\begin{array}{l}
aba = (c,a);\\
aca = (d,a);\\
ada = (b,1).\\
\end{array}
\right.
\end{equation}
Therefore, the restriction of $\psi$ to $St_\Gamma(1)$ gives a monomorphism
    $$\psi: St_\Gamma(1) \rightarrow \Gamma \times \Gamma,$$
which is not onto.
If $g \in St_\Gamma(1)$ is represented by a reduced word $w \in X^\ast$ then one can easily find the reduced forms of the automorphisms $g_0$ and $g_1$.
Indeed, in this case one can assume that $w \in \CR^e$ and
represented $w$ as a product
\begin{equation} \label{eq:st_product}
   w =   u_0 \cdot (a u_1 a) \cdot u_2 \cdot (a u_3 a) \cdot u_4 \ldots \cdot  \cdot u_{k-2} \cdot (a u_{k-1} a) \cdot u_k,
\end{equation}
where $u_0, \ldots, u_k \in \{b,c,d\}$ and $u_0$, $u_k$ are, perhaps, trivial. We refer to these $u_i$ and $(a u_j a)$  as to the {\em factors} of $w$. Now define two mappings $\phi_i: \CR^e \rightarrow \CR, i = 1,2,$ inductively on the number of factors. First, define $\phi_i$ on the factors according to the formulas (\ref{eq:st_split}):
\begin{equation}\label{eq:st_split-zero}
\left\{
\begin{array}{llllll}
\phi_0(b) = a&&&&&\\
\phi_0(c) = a\\
\phi_0(d) = \varepsilon\\
\end{array}
\right.
\left\{
\begin{array}{l}
\phi_0(aba) = c\\
\phi_0(aca) = d\\
\phi_0(ada) = b.\\
\end{array}
\right.
\end{equation}
\begin{equation}\label{eq:st_split-one}
\left\{
\begin{array}{llllll}
\phi_1(b) = c&&&&&\\
\phi_1(c) = d\\
\phi_1(d) = b\\
\end{array}
\right.
\left\{
\begin{array}{l}
\phi_1(aba) = a\\
\phi_1(aca) = a\\
\phi_1(ada) = 1.\\
\end{array}
\right.
\end{equation}
Then define by induction
    $$w_i = \phi_i(w) = red( \phi_i(v_1)\phi_i(v_2 \ldots v_k)), \ \ \ i = 0,1,$$
where $w = v_1 \ldots v_k$ is the factor decomposition (\ref{eq:st_product})  of $w$.
It follows immediately from the construction that for any $w \in \CR^e$
    $$ w ~\stackrel{\psi}{\mapsto}~ (w_0,w_1).$$
Notice that it takes time $O(|w|)$ to compute the pair $(w_0,w_1)$.

\begin{lemma}\label{le:split_length}
Let $w \in \CR^e$ and $w = (w_0,w_1)$. Then:
 \begin{itemize}
 \item [1)]   $$|w_0|,|w_1| \le \frac{|w|+1}{2}.$$
 \item [2)] Moreover, if $w$ starts with $a$ then $$|w_0|,|w_1| \le \frac{|w|}{2}.$$
\end{itemize}
\end{lemma}
\begin{proof}
Follows directly from the construction of $\phi_i, i = 1,2$ and the formulas (\ref{eq:st_split}) and (\ref{eq:st_split_prod}).
\end{proof}

\begin{remark}
\label{re:conjugation}
Let $w \in \CR^e$. Then conjugating, if necessary,  $w$ by its first letter or by its first two letters, and then reducing the result, one gets a word $w^\prime \in \CR^e$ which begins with $a$ and does not end on $a$.

\end{remark}

\subsection{The Word Problem in $\Gamma$}

Following  \cite{Harpe,Grig2006}, in this section  we discuss
an algorithm for the Word Problem in $\Gamma$. The algorithm is
based on three observations:
\begin{itemize}
    \item
If $|w|_a$ is odd then $w \not \in St_\Gamma(1)$, hence $w \neq 1$ in $\Gamma$.
    \item
If $|w|_a$ is even then  $w \in \CR^e$ and $w = (\phi_0(w),\phi_1(w))$. Moreover, since $\psi$ is a monomorphism, in this event
we have
    $$w=_\Gamma 1  ~~\Longleftrightarrow~~ \phi_0(w)=_\Gamma 1 ~\& ~ \phi_1(w)=_\Gamma 1$$
(here and below $w=_\Gamma 1$ means that $w = 1$ in $\Gamma$).
Therefore, the Word Problem for $w$ reduces to the Word Problem for $\phi_0(w)$
and $\phi_1(w)$, i.e., the process {\em splits}, or {\em branches}.
    \item
If $w = (w_0,w_1)$ then $|w_0|, |w_1| <  |w|$. Thus, the process stops in finitely many steps.
\end{itemize}
It is convenient to visualize the corresponding algorithm as an algorithm that
on an input $w \in X^\ast$ constructs a finite labeled rooted binary tree $T_w$.

\begin{algorithm}\label{al:WP}{\bf(Constructing the Decision Tree $T_w$)}
    \\{\sc Input.}
$w \in X^\ast$.
    \\{\sc Output.}
A finite labeled  rooted binary tree $T_w$.
    \\{\sc Computations.}
\begin{itemize}
    \item[A.]
{\bf(Initialization)}
Let  $T_0$ be a rooted binary tree with a single vertex (the root) $w$.
    \item[B.]
{\bf(Verification)} Let $T$ be a current rooted binary tree whose vertices are words in $(X\cup X^{-1})^\ast$ and some of them are marked by ``Yes'' or ``No''. Let $u$ be unmarked leaf in $T$. Then

Compute $|u|_a$. If $|u|_a$ is odd, then label $u$ by ``No'' and output the resulting tree as $T_w$.

Otherwise, compute $red(u)$ and take its conjugate $red(u)^\prime$ form Remark \ref{re:conjugation}. If $red(u)$ is empty then label $u$ by ``Yes'' and go to step B. Otherwise, go to C.

If there is no unmarked leaves in $T$ output $T$ as $T_w$.

    \item [C.]
{\bf(Splitting)}
Compute $\phi_0(u)$ and  $\phi_1(u)$ and add them as the ``left'' and the ``right'' children of $u$. Go to B.
\end{itemize}

 \end{algorithm}

\begin{proposition}
\label{pr:decision-tree-WP}
For a given word $w \in (X \cup X^{-1})^\ast$. The height of the tree $T_w$ is at most $\log_2|w| +1$.
\end{proposition}
 \begin{proof}
 By Lemma \ref{le:split_length} $|w_i| \leq |w|/2$, $i = 0,1$, for $w \in \CR^e$,  provided it begins with $a$ but does not end on $a$. Hence, starting with $w$ the Algorithm \ref{al:WP} can make at most $\log_2|w|$ splittings, since it does not split empty words. The verification step does not increase the height.
 \end{proof}

The following result is known (see, for example, \cite{Harpe,Grig2006}), but we need the proof for  the sake of references.
\begin{theorem} {\bf (Word Problem)}
The computational complexity of the Word Problem for the group $\Gamma$
is bounded by $O(n \log_2 n)$.
\end{theorem}
\begin{proof}
The algorithm for WP in $\Gamma$ works as follows. Given $w \in (X \cup X^{-1})^\ast$ it computes, first, the decision tree $T_w$. If $T_w$ has a vertex marked by ``No'' then $w \neq  1$ in $\Gamma$, otherwise $w = 1$ in $\Gamma$. By Proposition \ref{pr:decision-tree-WP} the decision tree $T_w$ has at most $\log_2|w|$ levels. Hence, to estimate the time required for the algorithm to construct $T_w$ one needs only to bound the time required to construct an arbitrary level in $T_w$. The verification step, as well as the splitting step, at a leaf $u$ requires only linear time in $|u|$. The total length of the vertices at a given level in $T_w$ is at most $|w|$. Hence the upper time bound for the complexity is $O(|w|\log_2|w|)$ as claimed.
\end{proof}

\subsection{The subgroup $K$}

As we mentioned above, the monomorphism $\psi:St_\Gamma(1) \rightarrow \Gamma \times \Gamma$ is not onto. In this section we describe a method how one can effectively check if a given pair  $(w_0,w_1) \in \Gamma \times \Gamma$ has a pre-image under $\psi$, and, if so, to compute it.  We refer to  \cite[Sections VIII.30, VIII.25]{Harpe} for details.

\begin{figure}[h]
\centerline{\includegraphics[scale=1]{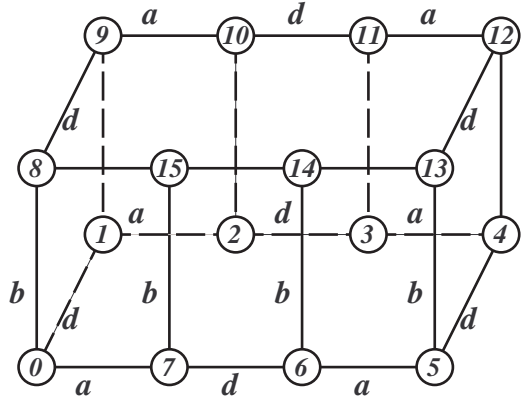} }
\caption{Schreier graph of $K \le \Gamma$ relative to generators $\{a,b,d\}$.}
 \label{fi:Schreier_K}
\end{figure}

Let $K$ be  the normal subgroup  of $\Gamma$ generated by the element $abab$
    $$K = \gp{abab}^\Gamma.$$
It turns out that $K$ has index $16$ in $\Gamma$  and
    $$K = \gp{abab, badabada, abadabad}.$$
The Schreier coset graph  of $K$ is shown in Figure
\ref{fi:Schreier_K}.
We denote the coset representatives of $K$ in $\Gamma$ by $1=g_0,g_1, \ldots, g_{15}$
according to the numbers in Figure \ref{fi:Schreier_K}. Observe, that $K$ is a subgroup of $St_\Gamma(1)$ of index 8 with coset representatives $g_0, g_1, g_4, g_5, g_8, g_9, g_{12}, g_{13}$.

\begin{lemma}\label{le:K_split_component}
For any $k\in K$ there exist elements $u, v  \in K$ such that
$u = (k,1)$ and $v = (1,k)$. In particular, $\psi(K) \geq K \times K$.
\end{lemma}

\begin{proof}
It is sufficient to prove the statement for the generators $abab, badabada, abadabad$ of the subgroup $K$. A straightforward verification shows that
    $$b\cdot ada\cdot b\cdot ada = (abab,1),$$
    $$badab \cdot aca \cdot badab \cdot aca = (abadabad,1),$$
    $$c \cdot badab \cdot aca \cdot badab \cdot aca \cdot c = (badabada,1).$$
Observe, that the words $w$ in the left-hand sides of the equalities above represent elements from $K$. Indeed, starting at the vertex $0$ and reading such a word $w$ in the Schreier graph above (beforehand replacing $c$ with $bd$) one ends up again at $0$,  thus proving the claim.
\end{proof}

\begin{lemma}
\label{le:image-psi}
Let  $D = \gp{(1,d),(1,a)} \leq \Gamma \times \Gamma$. Then:
\begin{itemize}
\item  [1)] $D$ is isomorphic to the Dihedral group of order $8$.
\item [2)] $\Gamma \times \Gamma = \psi(St_\Gamma(1)) \rtimes D$.
\end{itemize}
\end{lemma}
\begin{proof}
See \cite{Harpe}, page 229.
\end{proof}

\begin{figure}[h]
\centerline{\includegraphics[scale=0.8]{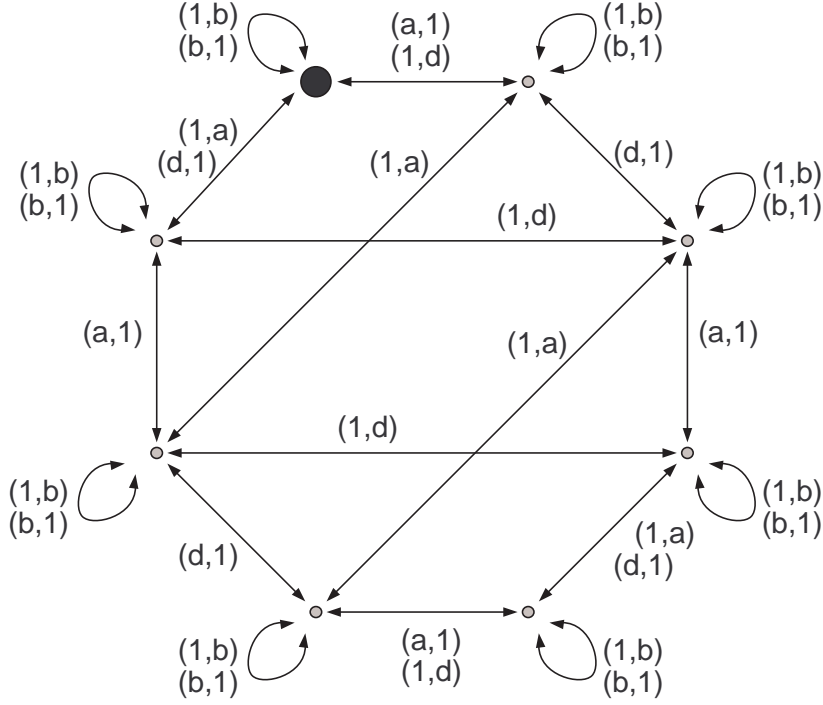} }
\caption{Schreier coset graph  of $\psi(St_\Gamma(1)) \le \Gamma \times \Gamma$ relative to the generating set
$\{(1,a),(1,b),(1,d), (a,1), (b,1), (d,1)\}$. The big black dot corresponds to the
coset $St_\Gamma(1)$.}
 \label{fi:Schreier-psi}
\end{figure}

The Schreier coset graph  of $\psi(St_\Gamma(1)) \le \Gamma \times \Gamma$ is shown in Figure
\ref{fi:Schreier-psi}.

\begin{lemma}\label{le:K_split_preimage}
Let $(u_0,u_1),  (v_0,v_1) \in \Gamma \times \Gamma$ be such that
$K u_0 = K v_0$ and $Ku_1=Kv_1$. If there exists $u\in\Gamma$ such that
$\psi(u) = (u_0,u_1)$, then there exists $v\in\Gamma$ such that
$\psi(v) = (v_0,v_1)$. Moreover, $Ku=Kv$.
\end{lemma}

\begin{proof}
Indeed, let $u_0 = k_0v_0, u_1 = k_1v_1$ for some $k_1,k_2 \in K$. Then
    $$(u_0,u_1) = (k_0v_0,k_1v_1) = (k_0,1)(1,k_1)(v_0,v_1).$$
By  Lemma \ref{le:K_split_component} $(k_0,1)$ and $(1,k_1)$ have pre-images
in $K$ under $\psi$. Therefore if $(u_0,u_1)$ has a pre-image in $St_\Gamma(1)$ then $(v_0,v_1)$ also has a pre-image in $St_\Gamma(1)$ and these pre-images lie in the same $K$-coset,
as required.
\end{proof}

Table \ref{tb:split} below describes completely the $K$-cosets of the pre-images $w$ under $\psi$ of pairs $(w_0,w_1)$ of elements from $\Gamma \times \Gamma$ if the cosets of the components $w_0, w_1$ are given (here numbers $i$ are the indices of the representatives $g_i$ of the cosets of $K$).


  \begin{table}[ht]
  \centering
  \begin{tabular*}{1\textwidth}%
     {@{\extracolsep{\fill}}|l|c||l|c||l|c||l|c|}
  \hline
$(w_0,w_1)$ & $w$ & $(w_0,w_1)$ & $w$ & $(w_0,w_1)$ & $w$ & $(w_0,w_1)$ & $w$ \\
  \hline
$( 0, 0)$ & $ 0$ & $( 8, 0)$ & $ 5$ & $( 4, 4)$ & $ 0$ & $(12, 4)$ & $ 5$ \\
  \hline
$( 0, 8)$ & $ 1$ & $( 8, 8)$ & $ 4$ & $( 4,12)$ & $ 1$ & $(12,12)$ & $ 4$ \\
  \hline
$( 1, 7)$ & $13$ & $( 9, 7)$ & $ 8$ & $( 5, 3)$ & $13$ & $(13, 3)$ & $ 8$ \\
  \hline
$( 1,15)$ & $12$ & $( 9,15)$ & $ 9$ & $( 5,11)$ & $12$ & $(13,11)$ & $ 9$ \\
  \hline
$( 2, 6)$ & $ 4$ & $(10, 6)$ & $ 1$ & $( 6, 2)$ & $ 4$ & $(14, 2)$ & $ 1$ \\
  \hline
$( 2,14)$ & $ 5$ & $(10,14)$ & $ 0$ & $( 6,10)$ & $ 5$ & $(14,10)$ & $ 0$ \\
  \hline
$( 3, 5)$ & $ 9$ & $(11, 5)$ & $12$ & $( 7, 1)$ & $ 9$ & $(15, 1)$ & $12$ \\
  \hline
$( 3,13)$ & $ 8$ & $(11,13)$ & $13$ & $( 7, 9)$ & $ 8$ & $(15, 9)$ & $13$ \\
  \hline
  \end{tabular*}
  \caption{}
  \label{tb:split}
  \end{table}

\begin{lemma}\label{le:lift}
Let $v_0,v_1 \in \Gamma$, $K v_0 = K g_i$ and $K v_1 = K g_j$. Then there
exists $v\in St_\Gamma(1)$ such that $v = (v_0,v_1)$ if and only if the pair $(i,j)$
is listed in the Table \ref{tb:split}.
\end{lemma}

\begin{proof} By Lemma \ref{le:K_split_preimage} the answer to the question whether a
pair $(u_0,u_1)$ has a pre-image in $St_\Gamma(1)$ under $\psi$
depends only on the coset $(K \times K)(u_0,u_1)$, which is completely determined by the cosets $Ku_0, Ku_1$
of the components. Therefore, it suffices to check which of the pairs $(g_i,g_j), i,j = 0, \ldots , 15$ lie in the subgroup $\psi(St_\Gamma(1))$. This can be easily done using the Schreier coset graph for $\psi(St_\Gamma(1)$, see Figure \ref{fi:Schreier-psi}.

\end{proof}

\section{Splittings}

In this section for a word $w \in \CR^e$ we study the metric properties of the {\em splitting} $w \rightarrow (w_0,w_1)$, where $w_i = red(\phi_i(w)), i= 1,2$. Namely, following \cite{Bartholdi:1998}, we study relations between norms (i.e., weighted lengths) of $w, w_0, w_1$.

Recall that
\begin{equation}\label{eq:basic}
\left\{
\begin{array}{llllll}
b = (a,c)&&&&&\\
c = (a,d)\\
d = (1,b)\\
\end{array}
\right.
\left\{
\begin{array}{l}
aba = (c,a)\\
aca = (d,a)\\
ada = (b,1).\\
\end{array}
\right.
\end{equation}

Let $\gamma_a,\gamma_b,\gamma_c,\gamma_d$ be fixed positive real values, termed   {\em weights}.
For a word $w \in X^\ast$ the number
    $$||w|| = \gamma_a \delta_a(w) + \gamma_b \delta_b(w) + \gamma_c \delta_c(w) + \gamma_d \delta_d(w).$$
is called the {\em norm} of $w$. The length $|w|$ is a special case of the norm when $\gamma_a=\gamma_b=\gamma_c=\gamma_d=1$.
In the following lemma we gather together some simple properties of the norm $||\cdot||$.

\begin{lemma}\label{le:norm_length}
Let  $u,v,w \in X^\ast$. Then  the following hold:
     \begin{itemize}
     \item [1)] $\min\{\gamma_a,\gamma_b,\gamma_c,\gamma_d\} \cdot |w| \le ||w|| \le \max\{\gamma_a,\gamma_b,\gamma_c,\gamma_d\} \cdot |w|.$
     \item [2)]  $||uv|| = ||u|| + ||v||.$
     \item [3)] If the numbers $\gamma_b,\gamma_c,\gamma_d$ satisfy the triangular inequality  then
    $$||red(w)|| \le ||w||.$$
    \end{itemize}
\end{lemma}

\begin{proof}
Straightforward verification.
\end{proof}

We define the weights  $\gamma_a, \gamma_d, \gamma_c, \gamma_b$ which will be in  use for the rest of the paper.
Let $\al$ be the unique real root of the polynomial $2 x^3 - x^2 - x - 1$,
    $$\al \in (1.233751,1.233752) \approx 1.23375.$$
Put
\begin{align*}
\gamma_a &= \al^2 + \al - 1 \in (1.755892,1.755896) \approx 1.7559, \\
\gamma_b &= 2, \\
\gamma_c &= \al^2 - \al + 1 \in (1.28839,1.288392) \approx 1.288, \\
\gamma_d &= -\al^2 + \al + 1 \in (0.711608,0.71161) \approx 0.712.
\end{align*}
Obviously, the weights $\gamma_a, \gamma_d, \gamma_c, \gamma_b$ satisfy the triangle inequality. Notice that the weights used in \cite{Bartholdi:1998} are slightly different (though related).

\begin{lemma}
\label{le:weights-rel}
The  following equalities hold:
\begin{align}
||a|| + ||b|| &= \al(||a||+||c||), \notag \\
||a|| + ||c|| &= \al(||a||+||d||), \label{eq:weights-rel}\\
||a|| + ||d|| &= \al||b||.  \notag
\end{align}
 \end{lemma}
 \begin{proof}
 Straightforward verification.
 \end{proof}

\begin{remark}
The choice of the weights  $\gamma_a, \gamma_b, \gamma_c, \gamma_d$  is optimal in the following sense: the value of $\alpha$ is maximal  for the weights satisfying  the triangle inequality (for $\gamma_b, \gamma_c, \gamma_d$) and relations \eqref{eq:weights-rel}.
\end{remark}

 The following results  establish some relations between $||w||$ and $||w_0||, ||w_1||$.

\begin{lemma} \label{le:splitting-norm}
Let $w \in \CR^e$ and $w_i = red(\psi_i(w))$ for $i=0,1$. Then the following hold:

\begin{itemize}
 \item [1)] If $w$ is of the form $*a*a\dots*a$ or $a*a*\dots a*$ (where stars stand for letters $b,c,d$) then
$$
  \al (||w_0|| + ||w_1||) \le ||w||.
$$
  \item [2)] If $w$ is of the form $*a*a\dots*$ then
$$
  \al (||w_0|| + ||w_1||) \le ||w|| + ||a||.
$$
  \item [3)] If $w$ is of the form $a*a\dots*a$ then
$$
  \al (||w_0|| + ||w_1||) \le ||w|| - ||a||.
$$
\end{itemize}
\end{lemma}
\begin{proof}
To prove 1) suppose  $w$ is in the form $*a*a\dots*a$. If $w = * \cdot a*a$ then the routine case by case verification, based on \eqref{eq:basic} and  \eqref{eq:weights-rel},  shows that
$$\al (||w_0|| + ||w_1||) \le \al (||\psi_0(w)|| + ||\psi_1(w)||) =  ||w||.$$
  In general   $w$  can be presented as a product of factors of the type $* \cdot a*a$:
$$ w = (x_1 \cdot ay_1a) \cdot \ldots \cdot (x_k\cdot ay_ka),$$
where $x_i, y_i \in \{b,c,d\}$, $i= 1, \ldots,k$. In this case   Lemma \ref{le:norm_length}, item  2), gives (using the fact that $\psi_0, \psi_1$ are homomorphisms)
$$
  \al (||\psi_0(w)|| + ||\psi_1(w)||) =  \al(\sum_{i = 1}^k ||\psi_0(x_i\cdot ay_ia)|| + \sum_{i = 1}^k ||\psi_1(x_i\cdot ay_ia)||) = \sum_{i = 1}^k ||x_i\cdot ay_ia|| =  |w||.
$$
Hence, by Lemma \ref{le:norm_length}, item 3),  $\al (||w_0|| + ||w_1||) \le ||w||$, as claimed.

To show 2) observe first that for $w = *a*a*$  one has
$$\al (||w_0|| + ||w_1||) \le ||w|| +||a||.$$
 Now the result follows from this and the case 1) above.

To see 3) it suffices to notice that for $w = a*a$
$$
  \al (||\psi_0(w)|| + ||\psi_1(w)||) \leq ||w|| -||a||
$$
and then apply an argument as above.
\end{proof}

\begin{corollary} \label{co:splitting-norm}
Let $w \in \CR$. For $i = 0,1$ put
$$w_i =
\left\{
\begin{array}{ll}
red(\psi_i(w)), & \mbox{if } w \in \CR^e;\\
red(\psi_i(wa)), & \mbox{if } w \not\in \CR^e.\\
\end{array}
\right.
$$
Then
$$\al (||w_0|| + ||w_1||) \le ||w|| + ||a||.$$
\end{corollary}
  \begin{proof}
  If $w \in \CR^e$ then the result follows directly from Lemma \ref{le:splitting-norm}. Suppose that $w \not \in \CR^e$. There are two cases to consider.
     Case 1. $w$ ends on a letter from $\{b,c,d\}$, i.e., $ w= u\cdot *$. Then $wa = u\cdot * \cdot a$ and if $u$ starts with $*$ then by Lemma \ref{le:splitting-norm}, item 1),
     $$\al (||w_0|| + ||w_1||) \leq ||wa|| = ||w|| + ||a||,$$
     as required. Otherwise, $u$ starts with $a$, then by Lemma \ref{le:splitting-norm}, item 3),
     $$\al (||w_0|| + ||w_1||) \leq ||wa|| - ||a|| = ||w||,$$
  which implies the result.
   \end{proof}

\begin{proposition}\label{pr:split_growth}
Let $w \in \CR$.
For $i = 0,1$ put
$$w_i =
\left\{
\begin{array}{ll}
red(\psi_i(w)), & \mbox{if } w \in \CR^e;\\
red(\psi_i(wa)), & \mbox{if } w \not\in \CR^e.\\
\end{array}
\right.
$$
Then the following hold:
 \begin{itemize}
 \item [1)]  If $\|w\|\ge 9$  then  $$\frac{||w||}{||w_0||+||w_1||}  \ge 1.03.$$
   \item If $||w||\ge 200$ then
    $$\frac{||w||}{||w_0||+||w_1||}  \ge 1.22.$$
    \end{itemize}
\end{proposition}

\begin{proof}
By Corollary~\ref{co:splitting-norm}
$$\al (||w_0|| + ||w_1||) \le ||w|| + ||a||. $$
Hence, if $||a||<0.01 ||w||$, which is the case when $||w||\ge 200$, then  the second inequality holds.
Similarly,  if $||w||\ge 9$ then the first inequality holds.
\end{proof}

\section{The Conjugacy Problem in the Grigorchuk Group}

In this section we prove that the Conjugacy Problem  (CP) in $\Gamma$ has a polynomial time decision algorithm.

\begin{lemma}\label{le:conjugacy_st}
If $u, v \in X^\ast$ are conjugate in $\Gamma$ then
$u \in St_\Gamma(1) ~~\Leftrightarrow~~ v \in St_\Gamma(1).$
\end{lemma}

\begin{proof}
Follows from Lemma \ref{le:st_generators} and the definition of conjugate elements.

\end{proof}

The next lemma describes behavior of conjugation relative to the splittings $w \rightarrow (w_0,w_1)$.  Below we frequently use the same notation for a word from $X^\ast$ and the element it represent in $\Gamma$, since it is clear from the context which one is which.

\begin{lemma}\cite{Grig2006} \label{le:conjugacy_split}
Let $u, v, x \in X^\ast$. Then the following hold in $\Gamma$:
\begin{itemize}
    \item[(S1)]
If $u,v,x \in St_\Gamma(1)$, and $u = (u_0,u_1)$, $v = (v_0,v_1)$, $x = (x_0,x_1)$ then
$$
    u = x^{-1} v x ~~\Leftrightarrow~~
    \left\{
    \begin{array}{l}
    u_0 = x_0^{-1} v_0 x_0,\\
    u_1 = x_1^{-1} v_1 x_1.\\
    \end{array}
    \right.
$$
    \item[(S2)]
If $u,v,xa \in St_\Gamma(1)$, and $u = (u_0,u_1)$, $v = (v_0,v_1)$, $xa = (x_0,x_1)$ then
$$
    u = x^{-1} v x ~~\Leftrightarrow~~
    \left\{
    \begin{array}{l}
    u_0 = x_1^{-1} v_1 x_1,\\
    u_1 = x_0^{-1} v_0 x_0.\\
    \end{array}
    \right.
$$
    \item[(N1)]
If $ua,va,x \in St_\Gamma(1)$, and $ua = (u_0,u_1)$, $va = (v_0,v_1)$, $x = (x_0,x_1)$ then
$$
    u = x^{-1} v x ~~\Leftrightarrow~~
    \left\{
    \begin{array}{l}
    u_0 u_1 = x_0^{-1} v_0 v_1 x_0,\\
    x_1 = v_1 x_0 u_1^{-1}.\\
    \end{array}
    \right.
$$
    \item[(N2)]
If  $ua,va,xa \in St_\Gamma(1)$, and $ua = (u_0,u_1)$, $va = (v_0,v_1)$, $xa = (x_0,x_1)$ then
$$
    u = x^{-1} v x ~~\Leftrightarrow~~
    \left\{
    \begin{array}{l}
    u_1 u_0 = x_0^{-1} v_0 v_1 x_0,\\
    x_1 = v_1 x_0 u_0^{-1}.\\
    \end{array}
    \right.
$$
\end{itemize}
\end{lemma}

\begin{proof}
(S1) and (S2) immediately follow  from  (\ref{eq:st_split_prod}) and (\ref{eq:st_split_conj}).

To see (N1) observe first that $u = uaa, v = vaa$ therefore
$$u = x^{-1}vx \Leftrightarrow  (u_0,u_1)a = (x_0^{-1},x_1^{-1})(v_0,v_1)a(x_0,x_1) \Leftrightarrow $$
 $$(u_0,u_1) = (x_0^{-1},x_1^{-1})(v_0,v_1)a(x_0,x_1)a \Leftrightarrow  (u_0,u_1) = (x_0^{-1},x_1^{-1})(v_0,v_1)(x_1,x_0).$$
 This implies equalities
 \begin{equation}
 \label{eq:N1}
 u_0 = x_0^{-1}v_0x_1, \ \ u_1 = x_1^{-1}v_1x_0
 \end{equation}
 Multiplying the equalities (\ref{eq:N1}) one gets $u_0u_1 = x_0^{-1}v_0v_1x_0$  - the first equality in (N1). Now from this equality one gets $x_0^{-1}v_0 =u_0u_1x_0^{-1}v_1^{-1}$. Substituting this into the first equality  in (\ref{eq:N1}) gives (after the standard manipulations) the second equality $x_1 = v_1x_0u_1^{-1}$  of (N1), as required.
 A similar argument proves (N2).
\end{proof}

 For a pair
of elements $u,v \in \Gamma$ define  a set
    $$Q(u,v) = \{g_i \mid \exists x\in \Gamma \mbox{ s.t. } u = x^{-1} v x ~\mbox{ and }~ Kx=Kg_i\},$$
where $g_0,\ldots,g_{15}$ are $K$-coset representatives of $\Gamma$  chosen above.  Clearly, $u$ and $v$ are conjugate in $\Gamma$ if and only if $Q(u,v) \neq \emptyset$.

The following is a key  lemma  in the solution of the conjugacy problem in $\Gamma$ (see \cite{Grig2006}).
\begin{lemma}
 \label{le:Q-u-v}
  Let $u, v \in X^\ast$.  Then the following hold:
\begin{itemize}
    \item [1)]
If $u = (u_0,u_1), v=(v_0,v_1) \in St_\Gamma(1)$ then
    $$Q(u,v) = \psi^{-1}[Q(u_0,v_0) \times Q(u_1,v_1)] \cup \psi^{-1}[Q(u_1,v_0) \times Q(u_0,v_1) ].$$
Moreover, it takes constant time to compute $Q(u,v)$ if the sets $Q(u_0,v_0)$,  $Q(u_1,v_1)$, $Q(u_0,v_1)$, and $Q(u_1,v_0)$ are given.
    \item [2)]
If $ua = (u_0,u_1), va=(v_0,v_1) \in St_\Gamma(1)$ then:
    $$Q(u,v) = \psi^{-1}\{ (g_i,g_j) \mid g_i\in Q(u_0u_1,v_0v_1) \mbox{ and } Kg_j = K v_1 g_i u_1^{-1} \} \cup $$
    $$\cup \psi^{-1}\{ (g_i,g_j) \mid g_i\in Q(u_1u_0,v_0v_1) \mbox{ and } Kg_j = K v_1 g_i u_0^{-1} \}.$$
Moreover, it takes constant time to compute $Q(u,v)$ if the sets $Q(u_0u_1,v_0v_1)$ and $Q(u_1u_0,v_0v_1)$ are given.
\end{itemize}
\end{lemma}

\begin{proof}
1) follows directly from  Lemma \ref{le:conjugacy_split}, items (S1) and (S2), and Lemma  \ref{le:lift}.

Similarly, 2) follows from  Lemma \ref{le:conjugacy_split}, items (N1) and (N2), and Lemma  \ref{le:lift}.

\end{proof}

Lemma \ref{le:Q-u-v} suggests a Branching  Decision Algorithm for the CP  (abbreviated to BDAC) in  $\Gamma$. The main idea of BDAC is the following: to check whether two given words $u, v \in X^\ast$  are conjugate  or not in $\Gamma$   it suffices to  verify if the set $Q(u,v)$  is empty or not. Hence
the conjugacy problem for elements $u, v$ is reduced to computing the set $Q(u,v)$. Now, to compute $Q(u,v)$ we are going to compute first the sets $Q(u^\prime,v^\prime)$ for a finite set of pairs $(u^\prime,v^\prime)$ that occur in the  branching process. To see how the  the algorithm works let $(u^\prime,v^\prime)$ be a current pair that occurs in BDAC. There are four cases to consider:
\begin{itemize}
    \item
{\bf (BDAC0)}
If $|u'|,|v'|\le 1$ then we use precomputed sets $Q(u',v')$. We compute them later in
Lemmas \ref{le:Q-obvious}, \ref{le:Q-4-6}, \ref{le:Qaa}, \ref{le:Qbb}, and \ref{le:Qdd}
    \item
{\bf (BDAC1)} If one of $u^\prime, v^\prime$  is in $St_\Gamma(1)$ and the other is not then by Lemma \ref{le:conjugacy_st}
the set $Q(u^\prime,v^\prime)$ is empty, in which case we mark the current pair $(u^\prime,v^\prime)$ by $\emptyset$.
    \item
{\bf (BDAC2)}
If $u^\prime = (u^\prime_0,u^\prime_1), v = (v^\prime_0,v^\prime_1) \in St_\Gamma(1)$ then, by Lemma \ref{le:Q-u-v} case 1), computation of $Q(u^\prime,v^\prime)$ reduces to  computing the sets $Q(u^\prime_0,v^\prime_0)$,  $Q(u^\prime_1,v^\prime_1)$, $Q(u^\prime_0,v^\prime_1)$, $Q(u^\prime_1,v^\prime_0)$.
    \item
{\bf (BDAC3)}
If $u^\prime a = (u^\prime_0,u^\prime_1), v^\prime a = (v^\prime_0,v^\prime_1) \in St_\Gamma(1)$ then, by Lemma \ref{le:Q-u-v} case 2), computation of $Q(u^\prime,v^\prime)$   reduces   to computing  the sets  $Q(u^\prime_0u^\prime_1,v^\prime_0v^\prime_1)$ and  $Q(u^\prime_1u^\prime_0,v^\prime_0v^\prime_1)$.
\end{itemize}
Thus, the process either assigns $\emptyset$ to the current pair of words  or branches at the pair (with four or two branches, depending on the case at hands).  By Proposition \ref{pr:split_growth},
each branching results in pairs of words with smaller norm, so the process eventually terminates in finitely many steps.

\begin{figure}[h]
\centerline{\includegraphics[scale=0.8]{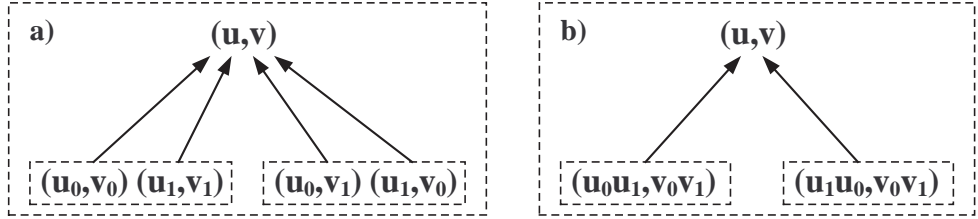} }
\caption{Two types of nodes in $T_{u,v}$ corresponding to cases BDAC2 and BDAC3.}
 \label{fi:nodes}
\end{figure}

To study the complexity of BDAC it is convenient to break it into two stages.

 {\bf Stage 1.} [Construction of the {\em Decision Conjugacy Tree} $T_{u,v}$.]
   At the first stage on an input $u, v \in X^\ast$ BDAC constructs  a finite labeled rooted tree $T_{u,v}$, where every vertex is a pair of words $(u_i,v_i)$  from $\CR$. Every vertex has degree at most four and some of them are  decorated with the symbol $\emptyset$. The pair $(u,v)$ is at the root of $T_{u,v}$. The construction of $T_{u,v}$ follows the rules BDAC0-BDAC3. Namely, if $(u^\prime,v^\prime)$ is a current node such that  $(u^\prime,v^\prime)$ falls into the case BDAC0 or BDAC1 then we leave this node as a leaf in $T_{u,v}$. If  $(u^\prime,v^\prime)$ falls into the case BDAC2 then the algorithm  constructs four children nodes $(u^\prime_0,v^\prime_0)$,  $(u^\prime_1,v^\prime_1)$, $(u^\prime_0,v^\prime_1)$ and  $(u^\prime_1,v^\prime_0)$ according to BDAC2 (see Figure \ref{fi:nodes} a)). If $(u^\prime,v^\prime)$ falls into the case BDAC3 then the algorithm constructs two  children nodes $(u^\prime_0u^\prime_1,v^\prime_0v^\prime_1)$ and  $(u^\prime_1u^\prime_0,v^\prime_0v^\prime_1)$ as described in BDAC3 (see Figure \ref{fi:nodes} b)).  Notice that it takes linear time in $|u| + |v|$ to assign $\emptyset$ to a given node or to produce its children.

 {\bf Stage 2.} [Computing the sets  $Q(u_i,v_i)$ involved in $T_{u,v}$.]
On the second stage BDAC, going from the leafs of  the tree $T_{u,v}$  to the root $(u,v)$, assigns to each vertex $(u^\prime,v^\prime)$
of $T_{u,v}$ the set $Q(u^\prime,v^\prime)$, computed as described in Lemma \ref{le:Q-u-v}.
It is left to assign the set $Q(u^\prime,v^\prime)$ to each leaf $(u^\prime,v^\prime)$ in $T_{u,v}$.
If the leaf $(u',v')$ corresponds to the case BDAC0 then we use precomputed sets $Q(u',v')$.
If the leaf $(u',v')$ corresponds to the case BDAC1 then $Q(u',v') = \emptyset$.
It takes constant time to assign the set $Q(u^\prime,v^\prime)$ to a vertex $(u^\prime,v^\prime)$, indeed, it is obvious
for the  leaves and follows from Lemma \ref{le:Q-u-v} for other vertices.

   It is clear from the description of BDAC that the time required for BDAC to stop and get the answer on an input $(u,v)$ can be estimated from above by the time to construct the tree $T_{u,v}$  and the time spent on  labeling  the vertices. Using the standard Breadth First algorithm the construction of the tree $T_{u,v}$ takes linear time in the size $|T_{u,v}|$ multiplied by the time spent at every vertex, so altogether is bounded from the above by $|T_{u,v}|(|u| + |v|)$.
   To get the polynomial estimate on the time complexity of BDAC we show below that the size of the tree $T_{u,v}$ is polynomial in terms of $|u| + |v|$.

The next result  shows  that for any words $u, v \in \CR$ the size $|T_{u,v}|$ of the tree $T_{u,v}$
 is polynomial in terms of $|u|$ and $|v|$   and gives estimates on the size.

\begin{lemma}\label{le:tree_size_9}
Let $u,v \in \CR$. If $\|v\| < 9$ and $\|u\| < 9$ then the size of the tree $T_{u,v}$ is not greater than $42$.
\end{lemma}

\begin{proof}
The set of pairs of words $(u,v)$ such that $\|v\| < 9$ and $\|u\| < 9$ is finite,
though relatively large, and, hence, the statement of the lemma can be checked by a straightforward verification.

Nevertheless, it is possible to check the correctness of the statement by hands.
Notice, that for any  child $(u',v')$ of $(u,v)$ the element $u'$ depends only on the element $u$
and the element $v'$ depends only on the element $v$. In other words, one  can assume  that
the left and the right words in the vertices are independent of each other, so one can consider each of them separately.
Table \ref{tb:T_size} contains all $95$ words $w$ of norm up to $9$.
For each $w$ it specifies the children of $w$ (defined as in cases BDAC2 and BDAC3)
and the size of the tree necessary to reach
words of lengths up to $1$ at the leaves. The greatest size is $21$ and hence, since we have
$2$ coordinates, the size of $T_{u,v}$ does not exceed $2\cdot 21 = 42$.

\end{proof}

\begin{proposition}
Let $u, v \in \CR$. Then the size of the tree $T_{u,v}$ is bounded by $2^{366} \rb{\max\{2|u|,2|v|\} }^{7}$.
\end{proposition}

\begin{proof}
Let $u, v \in \CR$. As in the proof of Lemma \ref{le:tree_size_9}  one   can consider the left and the right words in the vertices of the tree  separately.

By Lemma \ref{le:norm_length}, item 1), $\|w\|\le 2|w|$.
Consider an arbitrary branch
    $$(u_0,v_0), \ldots , (u_k,v_k)$$
in the tree $T_{u,v}$, starting at the root $(u_0,v_0) = (u,v)$ and ending at a leaf $(u_k,v_k)$.
It follows from Proposition \ref{pr:split_growth} that
for some $s \le \log_{1.22} \max\{2|u|,2|v|\}$ the inequalities $\|v_s\|\le 200$
and $\|u_s\| \le 200$ hold, i.e., any branch reaches a pair of words of norm up to $200$
in at most $\log_{1.22} \max\{2|u|,2|v|\}$ steps. Furthermore, by Proposition \ref{pr:split_growth}
there exists a number $t \le \log_{1.03} 200 < 180$ such that $\|v_{s+t}\| <9$ and $\|u_{s+t}\| < 9$.
By Lemma \ref{le:tree_size_9} the size of the tree $T_{v_{s+t},u_{s+t}}$ is not greater than $42$.
Thus, since the degree in each node in $T_{u,v}$ is not greater than $4$, it follows that
the size of $T_{u,v}$ is bounded by
    $$|T_{u,v}| \le 42\cdot 4^{180+ \log_{1.22} \max\{2|u|,2|v|\}} $$
    $$\le 2^{366} \rb{\max\{2|u|,2|v|\} }^{\log_{1.22} 4} \le 2^{366} \rb{\max\{2|u|,2|v|\} }^{7},$$
 as claimed.

\begin{table}[ht]
  \centering
  \newcolumntype{K}{>{\scriptsize$\columncolor[gray]{0.8}}l<{$}}
  \newcolumntype{L}{>{\scriptsize$}l<{$}}
  \begin{tabular*}{1\textwidth}%
     {@{\extracolsep{\fill}}|K|L||L|L||K|L||L|L||K|L||L|L|}
  \hline
w & w_1 & w_2 & \# & w & w_1 & w_2 & \# & w & w_1 & w_2 & \# \\
  \hline
  \hline
ab & ca & ac & 15 & adadad & b & b & 3 & cadac & aba & 1 & 5 \\
aba & c & a & 3 & ba & ac & ca & 15 & cadaca & aba & aba & 7 \\
abab & ca & ac & 15 & bab & 1 & 1 & 3 & cadad & ab & c & 17 \\
abac & ca & ad & 11 & baba & ac & ca & 15 & cadada & ad & cab & 9 \\
abaca & b & aba & 5 & babac & aca & cad & 21 & cadadad & ac & ca & 15 \\
abad & c & ab & 17 & babad & ac & cab & 13 & da & b & b & 3 \\
abada & cab & ad & 9 & bac & aba & b & 5 & dab & da & bac & 9 \\
abadad & dab & ac & 17 & baca & ad & ca & 11 & daba & c & ba & 17 \\
ac & da & ad & 7 & bacab & ada & cac & 7 & dabab & ca & bac & 13 \\
aca & d & a & 3 & bacac & ada & cad & 21 & dabac & ca & bad & 17 \\
acab & da & ac & 11 & bacad & ad & cab & 9 & dabad & c & bab & 5 \\
acaba & b & aba & 5 & bad & ad & cab & 9 & dabada & dab & bad & 19 \\
acac & da & ad & 7 & bada & ab & c & 17 & dac & ca & bad & 17 \\
acaca & 1 & 1 & 3 & badab & aba & 1 & 5 & daca & d & ba & 17 \\
acacad & dabad & b & 7 & badac & aba & b & 5 & dacab & da & bac & 9 \\
acad & d & ab & 17 & badad & ab & d & 17 & dacac & da & bad & 13 \\
acada & dab & ac & 17 & badada & ac & dab & 17 & dacaca & dabad & b & 7 \\
acadac & aba & aba & 7 & ca & ad & da & 7 & dacad & d & bab & 5 \\
acadad & cab & ad & 9 & cab & aba & b & 5 & dacada & cab & bac & 11 \\
ad & b & b & 3 & caba & ac & da & 11 & dacadad & dab & bad & 19 \\
ada & b & 1 & 3 & cabab & aca & dac & 21 & dad & 1 & 1 & 3 \\
adab & ba & c & 17 & cabac & aca & dad & 7 & dada & b & b & 3 \\
adaba & bac & da & 9 & cabad & ac & dab & 17 & dadab & ba & d & 17 \\
adabad & bad & dab & 19 & cac & 1 & 1 & 3 & dadaba & bad & ca & 17 \\
adac & ba & d & 17 & caca & ad & da & 7 & dadac & ba & c & 17 \\
adaca & bad & ca & 17 & cacab & ada & dac & 21 & dadaca & bac & da & 9 \\
adacac & b & dabad & 7 & cacac & ada & dad & 7 & dadacad & bad & dab & 19 \\
adacad & bac & cab & 11 & cacad & ad & dab & 13 & dadad & b & 1 & 3 \\
adad & b & b & 3 & cacada & b & dabad & 7 & dadada & b & b & 3 \\
adada & 1 & 1 & 3 & cad & ac & dab & 17 & dadadac & ca & ac & 15 \\
adadab & ca & bad & 17 & cada & ab & d & 17 & dadadad & 1 & 1 & 3 \\
adadac & da & bac & 9 & cadab & aba & b & 5 &&&& \\
  \hline
  \end{tabular*}
  \caption{}
  \label{tb:T_size}
  \end{table}

\end{proof}

\begin{theorem}
The Conjugacy Problem in the Grigorchuk group $\Gamma$ is decidable in $O(n^{8})$ time.
\end{theorem}

\begin{proof}
Given two words $u, v \in \CR$  the algorithm BDC construct the tree $T_{u,v}$, which size is bounded by
$2^{366} \rb{\max\{2|u|,2|v|\} }^{7}$.  Processing of each of the nodes of the tree  requires $O(\max\{|u|,|v|\})$ elementary steps.
Thus, the total complexity of the algorithm is bounded by $O\rb{\rb{\max\{|u|,|v|\} }^{8} }$.
\end{proof}

For completeness we list below the sets $Q(u,v)$ with $|u|, |v| \leq 1$.

 \begin{lemma}
 \label{le:Q-obvious}
 The following hold:
    $$Q(a,b) = Q(a,c) = Q(a,b) = \emptyset,$$
 $$Q(1,a) = Q(1,b) = Q(1,c) = Q(1,d) = \emptyset,$$
    $$Q(1,1) = \{0,\ldots, 15\}.$$
  \end{lemma}
\begin{proof}
Follows immediately from  Lemma \ref{le:conjugacy_st}.
\end{proof}

\begin{lemma} \label{le:Q-4-6}
    $$Q(b,c)=Q(b,d)=Q(c,d) = \emptyset.$$
\end{lemma}

\begin{proof}
Since, $b,c,d \in St_\Gamma(1)$ by Lemma \ref{le:Q-u-v} we have
    $$Q(c,d) = \psi^{-1}[Q(a,1) \times Q(d,b)] \cup \psi^{-1}[Q(a,b) \times Q(d,1)] = $$
    $$ = \psi^{-1}[\emptyset \times Q(d,b)] \cup \psi^{-1}[Q(a,b) \times \emptyset] = \emptyset.$$
Similarly,
    $$Q(b,d) = \psi^{-1}[Q(a,1) \times Q(c,b)] \cup \psi^{-1}[Q(a,b) \times Q(c,1)] = $$
    $$ = \psi^{-1}[\emptyset \times Q(d,b)] \cup \psi^{-1}[Q(a,b) \times \emptyset] = \emptyset.$$
Finally, using the obtained equalities we obtain
    $$Q(b,c) = \psi^{-1}[Q(a,a) \times Q(c,d)] \cup \psi^{-1}[Q(a,d) \times Q(c,a)] = $$
    $$ = \psi^{-1}[Q(a,a) \times \emptyset] \cup \psi^{-1}[\emptyset \times \emptyset] = \emptyset.$$
\end{proof}

\begin{lemma}\label{le:Qaa}
$Q(a,a) = \{0,3,4,7\}$.
\end{lemma}

\begin{proof}
Two cases to consider. If $x^{-1} a x = a$ and $x \in St_\Gamma(1)$ then $\psi(x)$ is of the form $(y,y) \in \Gamma\times \Gamma$. It
follows from Lemma \ref{le:lift} and Table \ref{tb:split} that in this case $Kx = Kg_0$
or $Kx = Kg_4$.

If $x^{-1} a x = a$ and $x \not \in St_\Gamma(1)$ then $\psi(xa)$ is of the form $(y,y)a \in \Gamma\times \Gamma$.
Using Figure \ref{fi:Schreier_K} it is easy to find that $K g_0a = Kg_7$ and $Kg_4a = Kg_3$.
Hence $Q(a,a) = \{0,3,4,7\}$.

\end{proof}

\begin{lemma} \label{le:Qbb}
$Q(b,b) = \{0,1,8,9\}$, $Q(c,c) = \{0,1,8,9\}$.
\end{lemma}

\begin{proof}
Notice that $1,b,c,d \in C_\Gamma(b)$, where $C_\Gamma(b)$ denotes the centralizer
of $b$ in $\Gamma$. Hence $Q(b,b) \supseteq \{0,1,8,9\}$.
On the other hand assume that $x^{-1} b x = b$. Then it is easy to see that $x=(x_0,x_1) \in St_\Gamma(1)$
and splitting $b$ into a pair $(a,c)$ we get
    $$Q(b,b) = \psi^{-1}[Q(a,a) \times Q(c,c)].$$
We proved in Lemma \ref{le:Qaa} that $Q(a,a) = \{0,3,4,7\}$.
Notice that in Table \ref{tb:split} all pairs $(i,j)$
with $i\in \{0,3,4,7\}$ define cosets with numbers $\{0,1,8,9\}$.
Hence $Q(b,b) \subseteq \{0,1,8,9\}$.

Proof for $Q(c,c)$ is the same.

\end{proof}

\begin{lemma} \label{le:Qdd}
$Q(d,d) = \{0,1,4,5,8,9,12,13\}$.
\end{lemma}

\begin{proof}
Notice that $1,b,c,d \in C_\Gamma(d)$, where $C_\Gamma(d)$ denotes the centralizer
of $d$ in $\Gamma$. Hence $Q(d,d) \supseteq \{0,1,8,9\}$
Furthermore, $ada,ada\cdot b,ada\cdot c,ada\cdot d \in C_\Gamma(d)$.
Hence $Q(d,d) \supseteq \{4,5,12,13\}$.

On the other hand assume that $x^{-1} d x = d$. Then it is easy to see that $x=(x_0,x_1) \in St_\Gamma(1)$
and splitting $d$ into a pair $(1,b)$ we get
    $$Q(d,d) = \psi^{-1} [Q(1,1) \times Q(b,b)].$$
Notice that in Table \ref{tb:split} all pairs $(i,j)$
with $i\in \{0,\ldots,15\}$ and $j \in \{0,1,8,9\}$ define cosets with numbers $\{0,1,4,5,8,9,12,13\}$.
Hence $Q(d,d) \subseteq \{0,1,4,5,8,9,12,13\}$.

\end{proof}

\end{document}